\documentclass[a4paper,12pt]{article}

 \usepackage{latexsym,amssymb,amsmath}

 \newtheorem{proposition}{Proposition}
 
 \newtheorem{definition}{Definition}

\begin{document}

\title{Monogenic Gaussian distribution in closed form and the Gaussian fundamental solution\thanks{This is a preprint of an article whose final and definitive form has been published in \textit{Complex Variables and Elliptic Equations}, 54 (2009), no. 5, 429--440.}}

\author{Dixan Pe\~na Pe\~na$^{\star,1}$ and Frank Sommen$^{\star\star,2}$}

\date{\normalsize{$^\star$Department of Mathematics, Aveiro University\\3810-193 Aveiro, Portugal\\
$^{\star\star}$Department of Mathematical Analysis, Ghent University\\9000 Ghent, Belgium}\\\vspace{0.4cm}
\small{$^1$e-mail: dixanpena@ua.pt; dixanpena@gmail.com\\
$^2$e-mail: fs@cage.ugent.be}} 

\maketitle

\begin{abstract}
\noindent In this paper we present a closed formula for the CK-extension of the Gaussian distribution in $\mathbb R^m$, and the monogenic version of the holomorphic function $\exp(z^2/2)/z$ which is a fundamental solution of the generalized Cauchy-Riemann operator.\vspace{0.2cm}\vspace{0.1cm}
\end{abstract}

\section{Introduction}

The main objects of study in Clifford analysis (see e.g. \cite{BDS,DSS,GuSp,K,KS,R}) are the so-called monogenic functions which may be described as null solutions of the Dirac operator, the latter being the higher dimensional analogue of the Cauchy-Riemann operator.

In this paper we deal with two very well-known techniques to generate monogenic functions: the Cauchy-Kowalevski extension (CK-extension) and Fueter's theorem. The first technique mentioned consists in monogenically extending analytic functions in $\mathbb R^m$ (see e.g. \cite{BDS,DSS,D,S,SJ}). The second one, named after the Swiss mathematician R. Fueter \cite{F}, gives a method to generate monogenic functions starting from a holomorphic function in the upper half of the complex plane (see \cite{KQS,LaLe,LaRa,D,DQS,DS,DS2,Q,QS,Sce,S3}).

The aim of this paper is to illustrate how Fueter's theorem may be used to derive two special functions in Clifford analysis: the monogenic Gaussian distribution in closed form, and the Gaussian fundamental solution which generalizes the complex fundamental solution $\exp(z^2/2)/z$.

This fundamental solution plays a key role in the theory of analytic functionals with unbounded carrier (see \cite{Mo}). The monogenic version is of the form 
\[\mathsf{E}(x)=\frac{\overline x}{\vert x\vert^{m+1}}+\mathsf{M}(x),\quad x\in\mathbb R^{m+1}\setminus\{0\},\]
$\mathsf{M}$ being an entire two-sided monogenic function. It satisfies an estimate of the form 
\[|\mathsf{E}(x)|\le C \exp\left(-|\underline x|^2/2\right),\quad|x_0|\le K,\quad|\underline x|\ge R,\] 
which is crucial for a monogenic generalization of the theory of analytic functionals with carrier in a strip domain.

\section{Clifford algebras and monogenic functions}

Clifford algebras were introduced in 1878 by the English geometer W.\ K.\ Clifford, generalizing the complex numbers and Hamilton's quaternions (see \cite{Cl}). They have important applications in geometry and theoretical physics.

We denote by $\mathbb R_{0,m}$ ($m\in\mathbb N$) the real Clifford algebra constructed over the orthonormal basis $(e_1,\ldots,e_m)$ of the Euclidean space $\mathbb R^m$. The basic axiom of this associative but non-commutative algebra is that the product of a vector with itself equals its squared length up to a minus sign, i.e. for any vector $\underline x=\sum_{j=1}^mx_je_j$ in $\mathbb R^m$, we have that
\[\underline x^2=-\vert\underline x\vert^2=-\sum_{j=1}^mx_j^2.\]
It thus follows that the elements of the basis submit to the multiplication rules
\begin{alignat*}{2}
e_j^2=-1&,&\qquad &j=1,\dots,m,\\
e_je_k+e_ke_j=0&,&\qquad &1\le j\neq k\le m.
\end{alignat*}
A basis for the algebra $\mathbb R_{0,m}$ is then given by the elements
\[e_A=e_{j_1}\cdots e_{j_k},\]
where $A=\{j_1,\dots,j_k\}\subset\{1,\dots,m\}$ is such that $j_1<\dots<j_k$. For the empty set $\emptyset$, we put $e_{\emptyset}=e_0=1$, the latter being the identity element. It follows that the dimension of $\mathbb R_{0,m}$ is $2^m$.

Any Clifford number $a\in\mathbb R_{0,m}$ may thus be written as
\[a=\sum_Aa_Ae_A,\quad a_A\in\mathbb R,\]
and its conjugate $\overline a$ is defined by
\[\overline a=\sum_Aa_A\overline e_A,\quad\overline e_A=(-1)^{\frac{k(k+1)}{2}}e_A,\quad\vert A\vert=k.\]
For each $k\in\{0,1,\dots,m\}$, we call
\[\mathbb R_{0,m}^{(k)}=\left\{a\in\mathbb R_{0,m}:\;a=\sum_{\vert A\vert=k}a_Ae_A\right\}\]
the subspace of $k$-vectors, i.e.\ the space spanned by the products of $k$ different basis vectors. In particular, the 0-vectors and 1-vectors are simply called scalars and vectors respectively.

Observe that $\mathbb R^{m+1}$ may be naturally embedded in the Clifford algebra $\mathbb R_{0,m}$ by associating to any element $(x_0,x_1,\ldots,x_m)\in\mathbb R^{m+1}$ the \lq\lq paravector" $x$ given by
\[x=x_0+\underline x.\]
Note that
\[\mathbb R_{0,m}=\bigoplus_{k=0}^m\mathbb R_{0,m}^{(k)}\]
and hence for any $a\in\mathbb R_{0,m}$
\[a=\sum_{k=0}^m[a]_k,\]
where $[a]_k$ is the projection of $a$ on $\mathbb R_{0,m}^{(k)}$. By means of the conjugation, a norm $\vert a\vert$ may be defined for each $a\in\mathbb R_{0,m}$ by putting
\[\vert a\vert^2=[a\overline a]_0=\sum_Aa_A^2.\]
Next, we introduce the Dirac operator
\[\partial_{\underline x}=\sum_{j=1}^me_j\partial_{x_j}\]
and the generalized Cauchy-Riemann operator
\[\partial_x=\partial_{x_0}+\partial_{\underline x}.\]
These operators factorize the Laplace operator in the sense that
\begin{equation}\label{fact1}
\Delta_{\underline x}=\sum_{j=1}^m\partial_{x_j}^2=-\partial_{\underline x}^2
\end{equation}
and
\begin{equation}\label{fact2}
\Delta_x=\partial_{x_0}^2+\Delta_{\underline x}=\partial_x\overline\partial_x=\overline\partial_x\partial_x.
\end{equation}

\begin{definition}
A function $f(\underline x)$ $($resp. $f(x)$$)$ defined and continuously differentiable in an open set $\Omega$ of $\,\mathbb R^m$ $($resp. $\mathbb R^{m+1}$$)$ and taking values in $\mathbb R_{0,m}$, is called a $($left$)$ monogenic function in $\Omega$ if and only if it fulfills in $\Omega$ the equation
\[\partial_{\underline x}f\equiv\sum_{j=1}^m\sum_Ae_je_A\partial_{x_j}f_A=0\quad(\text{resp.}\;\partial_xf\equiv\sum_{j=0}^m\sum_Ae_je_A\partial_{x_j}f_A=0).\]
\end{definition}
Note that in view of the non-commutativity of $\mathbb R_{0,m}$ a notion of right monogenicity may be defined in a similar way by letting act the Dirac operator or the generalized Cauchy-Riemann operator from the right. Functions which are both right and left monogenic are called two-sided monogenic functions.

Two basic examples of monogenic functions are $-\underline x/\vert\underline x\vert^m$ and $\overline x/\vert x\vert^{m+1}$, being these functions (up to a multiplicative constant) fundamental solutions of $\partial_{\underline x}$ and $\partial_x$, respectively. 

Finally, note that in view of (\ref{fact1}) and (\ref{fact2}) it follows that any monogenic function in $\Omega$ is harmonic in $\Omega$ and hence real-analytic in $\Omega$.

\section{Fueter's theorem}

Fueter's theorem was originally formulated in the setting of quaternionic analysis (see \cite{F}); and was later extended to the case of Clifford algebra-valued functions by Sce \cite{Sce}, Qian \cite{Q} and Sommen \cite{S3}. For other generalizations we refer the reader to \cite{KQS,LaLe,LaRa,D,DQS,DS,DS2,QS}.

Let $f(z)=u(x,y)+iv(x,y)$ ($z=x+iy$) be a holomorphic function in some open subset $\Xi\subset\mathbb C^+=\{z\in\mathbb C:\;y>0\}$; and let $P_k(\underline x)$ be a homogeneous monogenic polynomial of degree $k$ in $\mathbb R^m$, i.e.
\begin{alignat*}{2}
\partial_{\underline x}P_k(\underline x)&=0,&\quad\underline x&\in\mathbb R^m,\\
P_k(t\underline x)&=t^kP_k(\underline x),&\quad t&\in\mathbb R.
\end{alignat*}
Put $\underline\omega=\underline x/r$, with $r=\vert\underline x\vert$. In this paper, we are concerned with Sommen's generalization: \textit{if $m$ is an odd number, then the function}
\begin{equation*}\label{SoFueter}
\mathsf{Ft}\left[f(z),P_k(\underline x)\right](x)=\Delta_x^{k+\frac{m-1}{2}}\bigl[\bigl(u(x_0,r)+\underline\omega\,v(x_0,r)\bigr)P_k(\underline x)\bigr]
\end{equation*}
\textit{is monogenic in $\widetilde\Omega=\{x\in\mathbb R^{m+1}:\;(x_0,r)\in\Xi\}$}.

The proof of this generalization was based on the fact that
\[\bigl(u(x_0,r)+\underline\omega\,v(x_0,r)\bigr)P_k(\underline x)\]
may be written locally as $\overline\partial_x\big(h(x_0,r)P_k(\underline x)\big)$ for some $\mathbb R$-valued harmonic function $h$ of $x_0$ and $r$. Therefore using (\ref{fact2}), $\mathsf{Ft}\left[f(z),P_k(\underline x)\right]$ is monogenic if and only if (see \cite{S3})
\[\Delta_x^{k+\frac{m+1}{2}}\big(h(x_0,r)P_k(\underline x)\big)=0.\]
We notice that this version of Fueter's theorem provides us with the axial monogenic functions of degree $k$, i.e.
\[\mathsf{Ft}\left[f(z),P_k(\underline x)\right](x)=\bigl(A(x_0,r)+\underline\omega\,B(x_0,r)\bigr)P_k(\underline x),\]
where $A$ and $B$ are $\mathbb R$-valued and continuously differentiable functions in the variables $x_0$ and $r$ (see \cite{S4,S2}). It is not difficult to show that functions of this form are monogenic if and only if $A$ and $B$ satisfy the Vekua-type system
\begin{equation}\label{Veq}
\left\{\begin{array}{ll}\partial_{x_0}A-\partial_rB&=\displaystyle{\frac{2k+m-1}{r}}\,B\\\partial_{x_0}B+\partial_rA&=0.\end{array}\right.
\end{equation}
Using this fact, we presented in \cite{DQS} (see also \cite{D}) an alternative proof which has the advantage of allowing to compute some examples. Let us give an outline of the proof. We first showed by induction that
\[A=(2k+m-1)!!\,D_r\left(k+\frac{m-1}{2}\right)\{u\},\]
\[B=(2k+m-1)!!\,D^r\left(k+\frac{m-1}{2}\right)\{v\},\]
where $D_r(n)$ and $D^r(n)$ ($n\in\mathbb N_0$) are differential operators defined by
\begin{align*}
D_r(n)\{f\}&=\left(\frac{1}{r}\,\partial_r\right)^n\{f\},\quad\quad\quad\;\;\,D_r(0)\{f\}=f,\\
D^r(n)\{f\}&=\partial_r\left(\frac{D^r(n-1)\{f\}}{r}\right),\;D^r(0)\{f\}=f.
\end{align*}
Here we list some useful properties of these operators:
\begin{itemize}
\item[{\rm(i)}] $D^r(n)\{\partial_rf\}=\partial_rD_r(n)\{f\}$,
\item[{\rm(ii)}] $D_r(n)\{\partial_rf\}-\partial_rD^r(n)\{f\}=2n/r\,D^r(n)\{f\}$,
\item[{\rm(iii)}] $D_r(n)\{fg\}=\sum_{\nu=0}^n\binom{n}{\nu}D_r(n-\nu)\{f\}D_r(\nu)\{g\}$,
\item[{\rm(iv)}] $D^r(n)\{fg\}=\sum_{\nu=0}^n\binom{n}{\nu}D_r(n-\nu)\{f\}D^r(\nu)\{g\}$.
\end{itemize}
The final task was to prove that $A$ and $B$ satisfy the Vekua-type system (\ref{Veq}). In order to do that, it is necessary to use the assumptions on $u$ and $v$ and statements (i)-(ii).

Indeed,
\[\begin{split}
\partial_{x_0}&A-\partial_rB\\
&=(2k+m-1)!!\left(D_r\left(k+\frac{m-1}{2}\right)\{\partial_{x_0}u\}-\partial_rD^r\left(k+\frac{m-1}{2}\right)\{v\}\right)\\
&=(2k+m-1)!!\left(D_r\left(k+\frac{m-1}{2}\right)\{\partial_rv\}-\partial_rD^r\left(k+\frac{m-1}{2}\right)\{v\}\right)\\
&=\frac{2k+m-1}{r}\,(2k+m-1)!!\,D^r\left(k+\frac{m-1}{2}\right)\{v\}\\
&=\frac{2k+m-1}{r}\,B
\end{split}\]
and
\[\begin{split}
\partial_{x_0}&B+\partial_rA\\
&=(2k+m-1)!!\left(D^r\left(k+\frac{m-1}{2}\right)\{\partial_{x_0}v\}+\partial_rD_r\left(k+\frac{m-1}{2}\right)\{u\}\right)\\
&=(2k+m-1)!!\left(D^r\left(k+\frac{m-1}{2}\right)\{\partial_{x_0}v\}+D^r\left(k+\frac{m-1}{2}\right)\{\partial_ru\}\right)\\
&=(2k+m-1)!!\,D^r\left(k+\frac{m-1}{2}\right)\{\partial_{x_0}v+\partial_ru\}\\
&=0,
\end{split}\]
which completes the proof.

It is a simple matter to check that
\begin{alignat*}{2}
\mathsf{Ft}\left[cf(z),P_k(\underline x)\right]&=c\,\mathsf{Ft}\left[f(z),P_k(\underline x)\right],\quad c\in\mathbb R,\\
\mathsf{Ft}\left[f(z)+g(z),P_k(\underline x)\right]&=\mathsf{Ft}\left[f(z),P_k(\underline x)\right]+\mathsf{Ft}\left[g(z),P_k(\underline x)\right],
\end{alignat*}
where $f(z)$ and $g(z)$ are two holomorphic functions in the upper half of the complex plane. At this point it is important
to notice that $\mathsf{Ft}\left[if(z),P_k(\underline x)\right]\ne i\mathsf{Ft}\left[f(z),P_k(\underline x)\right]$. 

It is also worth remarking that for $k=0$, Fueter's theorem generates two-sided monogenic functions. More precisely,
\[\partial_x\mathsf{Ft}\left[f(z),1\right]=\mathsf{Ft}\left[f(z),1\right]\partial_x=0\;\;\text{in}\;\;\widetilde\Omega.\]
We will now compute some examples (see also \cite{DS2}).\vspace{0.1cm}

\noindent \textbf{Example 1.} Let $f(z)=iz=-y+ix$. It easily follows that
\begin{equation}\label{e1}
D_r(n)\{r\}=(-1)^{n+1}\frac{(2n-3)!!}{r^{2n-1}},
\end{equation}
\[D^r(n)\{x_0\}=(-1)^n\frac{(2n-1)!!}{r^{2n}}\,x_0.\]
We thus get the monogenic function
\begin{multline*}
\mathsf{Ft}\left[iz,P_k(\underline x)\right](x)=(-1)^{k+\frac{m-1}{2}}(2k+m-1)!!(2k+m-4)!!\\
\times\left(\frac{1}{r^{2k+m-2}}+\frac{(2k+m-2)x_0\underline x}{r^{2k+m}}\right)P_k(\underline x),\;\;\underline x\ne0.
\end{multline*}

\noindent \textbf{Example 2.} Consider
\[f(z)=\frac{1}{z}=\frac{x}{x^2+y^2}-i\frac{y}{x^2+y^2}.\]
It is easy to check that
\begin{equation}\label{e2}
D_r(n)\left\{\frac{x_0}{x_0^2+r^2}\right\}=(-1)^n\frac{2^nn!x_0}{(x_0^2+r^2)^{n+1}},
\end{equation}
\begin{equation}\label{e3}
D^r(n)\left\{\frac{r}{x_0^2+r^2}\right\}=(-1)^n\frac{2^nn!r}{(x_0^2+r^2)^{n+1}}.
\end{equation}
With this choice of initial function, we obtain the well-known monogenic function in $\mathbb R^{m+1}\setminus\{\underline x\ne0\}$:
\[\mathsf{Ft}\left[1/z,P_k(\underline x)\right](x)=(-1)^{k+\frac{m-1}{2}}((2k+m-1)!!)^2\left(\frac{\overline x}{\vert x\vert^{2k+m+1}}\right)P_k(\underline x).\]
Before introducing the last example, we first need to introduce the CK-extension technique. 

Let $f(\underline x)$ be an analytic function in $\mathbb R^m$. The CK-extension of $f$ is the unique monogenic extension $\mathsf{CK}[f]$ of $f$ to $\mathbb R^{m+1}$ and given by
\begin{equation*}\label{CK}
\mathsf{CK}[f](x)=\sum_{n=0}^\infty\frac{(-x_0)^n}{n!}\,\partial_{\underline x}^nf(\underline x).
\end{equation*}

\noindent \textbf{Example 3.} Let $n\in\mathbb N$ and take 
\begin{multline*}
f(z)=z^n\\
=\sum_{\nu=0}^{[n/2]}(-1)^\nu\binom{n}{2\nu}x^{n-2\nu}y^{2\nu}+i\sum_{\nu=0}^{[(n-1)/2]}(-1)^\nu\binom{n}{2\nu+1}x^{n-(2\nu+1)}y^{2\nu+1}.
\end{multline*}
For this initial function we have that 
\[u(x_0,r)+\underline\omega\,v(x_0,r)=\sum_{\nu=0}^n\binom{n}{\nu}x_0^{n-\nu}\underline x^\nu.\]
Therefore $\mathsf{Ft}\left[z^n,P_k(\underline x)\right]$ is a homogeneous monogenic polynomial of degree $n-k-m+1$ in $\mathbb R^{m+1}$. Moreover, 
\[\mathsf{Ft}\left[z^n,P_k(\underline x)\right](x)\Big\vert_{\underline x_0=0}=c\,\underline x^{n-(2k+m-1)}P_k(\underline x),\quad c\in\mathbb R.\]
From the above we can claim that $\mathsf{Ft}\left[z^n,P_k(\underline x)\right]$ equals (up to a multiplicative constant) the CK-extension of $\underline x^{n-(2k+m-1)}P_k(\underline x)$.

\begin{proposition}\label{prop1}
Let $H(z)=\sum_{n=0}^\infty c_nz^n$ $(c_n\in\mathbb R)$ be an entire function. Then $\mathsf{Ft}\left[H(z),P_k(\underline x)\right]$ is an entire monogenic function of the form
\[\sum_{n=2k+m-1}^\infty C_n\mathsf{CK}[\underline x^{n-(2k+m-1)}P_k(\underline x)](x),\quad C_n\in\mathbb R.\]
\end{proposition}

\section{Monogenic Gaussian distribution}

In this section, we will focus on the CK-extension of the Gaussian distribution $\exp(-\vert\underline x\vert^2/2)$ in $\mathbb R^m$. It may be given by the following series (see \cite{DSS})
\[\mathsf{CK}[\exp(-\vert\underline x\vert^2/2)](x)=\exp(-\vert\underline x\vert^2/2)\sum_{n=0}^\infty\frac{x_0^n}{n!}\,H_n(\underline x),\]
where the functions $H_n(\underline x)$ are polynomials in $\underline x$ of degree $n$ with real coefficients and satisfy the recurrence formula
\[H_{n+1}(\underline x)=\underline xH_n(\underline x)-\partial_{\underline x}H_n(\underline x).\]
The polynomials $H_n(\underline x)$ generalize the classical Hermite polynomials. It easily follows by induction that
\begin{align*}
H_{2n}(\underline x)&=\sum_{\nu=0}^n\binom{n}{\nu}c_n(\nu)\,\underline x^{2(n-\nu)},\\
H_{2n+1}(\underline x)&=\sum_{\nu=0}^n\binom{n}{\nu}c_{n+1}(\nu)\,\underline x^{2(n-\nu)+1},
\end{align*}
where
\[c_{n}(\nu)=\prod_{l=1}^\nu(m+2(n-l)),\quad c_{n}(0)=1.\]
From now on we will assume $m$ an odd number. It is clear that
\[\mathsf{CK}[\exp(-\vert\underline x\vert^2/2)](x)\Big\vert_{\underline x=0}=\sum_{n=0}^\infty\frac{x_0^{2n}}{(2n)!}\,c_{n}(n).\]
This series is the Taylor expansion of the function
\[\exp(x_0^2/2)\left(1+\sum_{n=1}^{\frac{m-1}{2}}\prod_{\nu=1}^n(m-(2\nu-1))\frac{x_0^{2n}}{(2n)!}\right).\]
Thus
\[\mathsf{CK}[\exp(-\vert\underline x\vert^2/2)](x)\Big\vert_{\underline x=0}=\exp(x_0^2/2)\left(1+\sum_{n=1}^{\frac{m-1}{2}}\prod_{\nu=1}^n(m-(2\nu-1))\frac{x_0^{2n}}{(2n)!}\right).\]
Consider the holomorphic function
\[f(z)=\exp(z^2/2)=\exp\left(\frac{x^2-y^2}{2}\right)\left(\cos(xy)+i\sin(xy)\right).\]
Let us now compute using Fueter's technique the corresponding monogenic function. It may be proved by induction that
\begin{align}
&D_r(n)\left\{\exp\left(\frac{x_0^2-r^2}{2}\right)\right\}=(-1)^n\exp\left(\frac{x_0^2-r^2}{2}\right)\label{e4},\\
&D_r(n)\{\cos(x_0r)\}=\sum_{\nu=1}^na_\nu^{(n)}\frac{x_0^\nu}{r^{2n-\nu}}\cos(x_0r+\nu\pi/2)\label{e5},\\
&D_r(n)\{\sin(x_0r)\}=\sum_{\nu=1}^na_\nu^{(n)}\frac{x_0^\nu}{r^{2n-\nu}}\sin(x_0r+\nu\pi/2)\label{e6},\\
&D^r(n)\{\sin(x_0r)\}=\sum_{\nu=0}^na_{\nu+1}^{(n+1)}\frac{x_0^\nu}{r^{2n-\nu}}\sin(x_0r+\nu\pi/2)\label{e7},
\end{align}
with 
\begin{align*}
a_1^{(n)}&=(-1)^{n+1}(2n-3)!!,\\
a_\nu^{(n+1)}&=-(2n-\nu)a_\nu^{(n)}+a_{\nu-1}^{(n)},\quad \nu=2,\dots,n,\\
a_n^{(n)}&=1.
\end{align*}
By statements (iii) and (iv), we see that
\begin{multline*}
D_r(n)\left\{\exp\left(\frac{x_0^2-r^2}{2}\right)\cos(x_0r)\right\}\\
=\exp\left(\frac{x_0^2-r^2}{2}\right)\sum_{\nu=0}^n\binom{n}{\nu}(-1)^{n-\nu}D_r(\nu)\{\cos(x_0r)\},
\end{multline*}
\begin{multline*}
D^r(n)\left\{\exp\left(\frac{x_0^2-r^2}{2}\right)\sin(x_0r)\right\}\\
=\exp\left(\frac{x_0^2-r^2}{2}\right)\sum_{\nu=0}^n\binom{n}{\nu}(-1)^{n-\nu}D^r(\nu)\{\sin(x_0r)\}.
\end{multline*}
Hence
\begin{multline*}\label{GaussDCA}
\mathsf{Ft}\left[\exp(z^2/2),P_k(\underline x)\right](x)=(2k+m-1)!!\\
\times\exp\left(\frac{x_0^2-r^2}{2}\right)\left(\sum_{\nu=0}^{k+\frac{m-1}{2}}\binom{k+\frac{m-1}{2}}{\nu}(-1)^{k+\frac{m-1}{2}-\nu}D_r(\nu)\{\cos(x_0r)\}\right.\\
\left.+\underline\omega\sum_{\nu=0}^{k+\frac{m-1}{2}}\binom{k+\frac{m-1}{2}}{\nu}(-1)^{k+\frac{m-1}{2}-\nu}D^r(\nu)\{\sin(x_0r)\}\right)P_k(\underline x),\quad\underline x\ne0.
\end{multline*}
Note that for $k=0$ the restriction of $\mathsf{Ft}\left[\exp(z^2/2),P_k(\underline x)\right]$ to $x_0=0$ is
\[(-1)^{\frac{m-1}{2}}(2k+m-1)!!\exp(-\vert\underline x\vert^2/2).\]
Therefore, for $k=0$, $\mathsf{Ft}\left[\exp(z^2/2),P_k(\underline x)\right]$ equals (up to a multiplicative constant) the CK-extension of $\exp(-\vert\underline x\vert^2/2)$ when $\underline x\ne0$.

\begin{proposition}
Let $m$ be an odd number. A closed formula for the CK-extension of the Gaussian distribution in $\mathbb R^m$ is given by
\begin{multline*}
\mathsf{CK}[\exp(-\vert\underline x\vert^2/2)](x)\\
=\left\{\begin{array}{ll}(-1)^{\frac{m-1}{2}}\displaystyle{\frac{\mathsf{Ft}\left[\exp(z^2/2),1)\right](x)}{(m-1)!!}}&\text{for}\quad\underline x\ne0,\\\displaystyle{\exp(x_0^2/2)\left(1+\sum_{n=1}^{\frac{m-1}{2}}\prod_{\nu=1}^n(m-(2\nu-1))\frac{x_0^{2n}}{(2n)!}\right)}&\text{for}\quad\underline x=0.\end{array}\right.
\end{multline*}
\end{proposition}
For the particular case $m=3$, we have that
\begin{multline*}
\mathsf{CK}[\exp(-\vert\underline x\vert^2/2)](x)=\exp\left(\frac{x_0^2-r^2}{2}\right)\bigg(\cos(x_0r)+\frac{x_0}{r}\sin(x_0r)\\
+\underline\omega\left(\sin(x_0r)+\frac{\sin(x_0r)}{r^2}-\frac{x_0}{r}\cos(x_0r)\right)\bigg),\;\text{for}\quad\underline x\ne0,
\end{multline*}
and
\[\mathsf{CK}[\exp(-\vert\underline x\vert^2/2)](x)=\exp(x_0^2/2)(1+x_0^2),\;\text{for}\quad\underline x=0.\]

\section{The Gaussian fundamental solution}

In this last section, we will test Fueter's technique with the initial holomorphic function
\begin{multline*}
f(z)=\frac{\exp(z^2/2)}{z}\\
=\frac{\exp\left(\frac{x^2-y^2}{2}\right)}{x^2+y^2}\Big(\big(x\cos(xy)+y\sin(xy)\big)+i\big(x\sin(xy)-y\cos(xy)\big)\Big).
\end{multline*}
Note that $\exp(z^2/2)/z$ may be written as $1/z+H(z)$, where $H(z)$ is a entire function whose coefficients in the Taylor expansion around $z=0$ are real. Thus, by Proposition \ref{prop1} and Example 2, there exists an entire two-sided monogenic function $\mathsf{M}$ such that
\[\mathsf{Ft}\left[\exp(z^2/2)/z,1\right](x)=(-1)^{\frac{m-1}{2}}((m-1)!!)^2\left(\frac{\overline x}{\vert x\vert^{m+1}}+\mathsf{M}(x)\right).\]
We note that $\mathsf{Ft}\left[\exp(z^2/2)/z,1\right]$ is explicitly given by 
\begin{multline*}
\mathsf{Ft}\left[\exp(z^2/2)/z,1\right](x)\\
=(m-1)!!\left(D_r\left(\frac{m-1}{2}\right)\left\{\frac{\exp\left(\frac{x_0^2-r^2}{2}\right)}{x_0^2+r^2}\big(x_0\cos(x_0r)+r\sin(x_0r)\big)\right\}\right.\\
\left.+\underline\omega\,D^r\left(\frac{m-1}{2}\right)\left\{\frac{\exp\left(\frac{x_0^2-r^2}{2}\right)}{x_0^2+r^2}\big(x_0\sin(x_0r)-r\cos(x_0r)\big)\right\}\right).
\end{multline*}
In view of the above and using statements (iii) and (iv), we can also assert that $\mathsf{Ft}\left[\exp(z^2/2)/z,1\right](x)$ may be expressed in terms of (\ref{e1})-(\ref{e7}) with $0\le n\le(m-1)/2$. Therefore
\[\mathsf{Ft}\left[\exp(z^2/2)/z,1\right](x)=\exp\left(\frac{x_0^2-r^2}{2}\right)\big(\alpha(x_0,r)+\underline\omega\,\beta(x_0,r)\big),\]
where $\alpha$ and $\beta$ are $\mathbb R$-valued functions.

Let $K,R>0$. We then have that for $\vert x_0\vert\le K$ and for $r\ge R$ the following inequality holds
\[\big\vert\alpha(x_0,r)+\underline\omega\,\beta(x_0,r)\big\vert\le C,\]
where $C$ denotes a positive constant depending on $K$ and $R$.

The above observations are summarized in the proposition below.

\begin{proposition}
Let $m$ be an odd number. The function
\[\mathsf{Ft}\left[\exp(z^2/2)/z,1\right](x)\]
equals (up to a multiplicative constant)
\[\frac{\overline x}{\vert x\vert^{m+1}}+\mathsf{M}(x),\]
where $\mathsf{M}(x)$ is an entire two-sided monogenic function. Moreover, for $\vert x_0\vert\le K$ and $r\ge R$, we have that
\[\big\vert\mathsf{Ft}\left[\exp(z^2/2)/z,1\right](x)\big\vert\le C \exp\left(-r^2/2\right).\]
\end{proposition}

\subsection*{Acknowledgments}

The first author was supported by a Post-Doctoral Grant of Funda\c{c}\~ao para a Ci\^encia e a Tecnologia (FCT), Portugal.


\begin{thebibliography}{99}

\bibitem[1]{BDS} F. Brackx, R. Delanghe and F. Sommen, \textit{Clifford analysis}, Research Notes in Mathematics, 76, Pitman (Advanced Publishing Program), Boston, MA, 1982.

\bibitem[2]{Cl} W.\ K.\ Clifford, \textit{Applications of Grassmann's Extensive Algebra}, Amer. J. Math. 1 (1878), no. 4, 350--358.

\bibitem[3]{DSS} R. Delanghe, F. Sommen and V. Sou\v cek, \textit{Clifford algebra and spinor-valued functions}, Mathematics and its Applications, 53, Kluwer Academic Publishers Group, Dordrecht, 1992.

\bibitem[4]{F} R. Fueter, \textit{Die funktionentheorie der differentialgleichungen $\Delta u=0$ und $\Delta\Delta u=0$ mit vier variablen}, Comm. Math. Helv. 7 (1935), 307--330.

\bibitem[5]{GuSp} K. G\"urlebeck and W. Spr\"ossig, \textit{Quaternionic and Clifford calculus for physicists and engineers}, Wiley and Sons Publ., 1997.

\bibitem[6]{KQS} K.\ I.\ Kou, T. Qian and F. Sommen, \textit{Generalizations of Fueter's theorem}, Methods Appl. Anal. 9 (2002), no. 2, 273--289.

\bibitem[7]{K} V.\ V.\ Kravchenko, \textit{Applied quaternionic analysis}, Research and Exposition in Mathematics, 28, Heldermann Verlag, Lemgo, 2003.

\bibitem[8]{KS} V.\ V.\ Kravchenko and M.\ V.\ Shapiro, \textit{Integral representations for spatial models of mathematical physics}, Pitman Research Notes in Mathematics Series, 351, Longman, Harlow, 1996.

\bibitem[9]{LaLe} G. Laville and E. Lehman, \textit{Analytic Cliffordian functions}, Ann. Acad. Sci. Fenn. Math. 29 (2004), no. 2, 251--268. 

\bibitem[10]{LaRa} G. Laville and I. Ramadanoff, \textit{Holomorphic Cliffordian functions}, Adv. Appl. Clifford Algebras 8 (1998), no. 2, 323--340.
    
\bibitem[11]{Mo} M. Morimoto, \textit{Analytic functionals with non-compact carrier}, Tokyo J. Math. 1 (1978), no. 1, 77--103.    

\bibitem[12]{D} D. Pe\~{n}a Pe\~{n}a, \textit{Cauchy-Kowalevski extensions, Fueter's theorems and boundary values of special systems in Clifford analysis}, Ph.D. Thesis, Ghent University, 2008.

\bibitem[13]{DQS} D. Pe\~{n}a Pe\~{n}a, T. Qian and F. Sommen, \textit{An alternative proof of Fueter's theorem}, Complex Var. Elliptic Equ. 51  (2006), no. 8-11, 913--922.

\bibitem[14]{DS} D. Pe\~{n}a Pe\~{n}a and F. Sommen, \textit{A generalization of Fueter's theorem}, Results Math. 49 (2006), no. 3-4, 301--311.

\bibitem[15]{DS2} D. Pe\~{n}a Pe\~{n}a and F. Sommen, \textit{A note on the Fueter theorem}, submitted for publication.

\bibitem[16]{Q} T. Qian, \textit{Generalization of Fueter's result to} $\mathbb R^{n+1}$, Atti Accad. Naz. Lincei Cl. Sci. Fis. Mat. Natur. Rend. Lincei (9) Mat. Appl. 8 (1997), no. 2, 111--117.

\bibitem[17]{QS} T. Qian and F. Sommen, \textit{Deriving harmonic functions in higher dimensional spaces}, Z. Anal. Anwendungen 22  (2003), no. 2, 275--288.

\bibitem[18]{R} J. Ryan, \textit{Basic Clifford analysis}, Cubo Mat. Educ. 2 (2000), 226--256.  

\bibitem[19]{Sce} M. Sce, \textit{Osservazioni sulle serie di potenze nei moduli quadratici}, Atti Accad. Naz. Lincei. Rend. Cl. Sci. Fis. Mat. Nat. (8) 23 (1957), 220--225.

\bibitem[20]{S4} F.  Sommen, \textit{Plane elliptic systems and monogenic functions in symmetric domains}, Rend. Circ. Mat. Palermo (2) 1984, no. 6, 259--269.

\bibitem[21]{S} F. Sommen, \textit{Monogenic functions on surfaces},  J. Reine Angew. Math. 361 (1985), 145--161.

\bibitem[22]{S2} F. Sommen, \textit{Special functions in Clifford analysis and axial symmetry}, J. Math. Anal. Appl. 130 (1988), no. 1, 110--133.

\bibitem[23]{S3} F. Sommen, \textit{On a generalization of Fueter's theorem}, Z. Anal. Anwendungen 19 (2000), no. 4, 899--902.

\bibitem[24]{SJ} F. Sommen and B. Jancewicz, \textit{Explicit solutions of the inhomogeneous Dirac equation}, J. Anal. Math. 71 (1997), 59--74.

\end{thebibliography}
\end{document}